\documentclass[12pt]{article}
\usepackage{mysty}
\usepackage{amsthm,amsmath,amssymb}
\usepackage{rotating}

\usepackage[colorlinks=true,citecolor=black,linkcolor=black,urlcolor=blue]{hyperref}

\newcommand{\arxiv}[1]{\href{http://arxiv.org/abs/#1}{\texttt{arXiv:#1}}}

\newcommand{\ve}[1]{\textit{\textbf{#1}}}
\newcommand{\abs}[1]{\lvert#1\rvert}
\newcommand{\sign}[1]{\textnormal{sign}(#1)}
\renewcommand{\H}[2]{\textnormal{H}_{#1}(#2)}
\newcommand{\R}{\mathbb R}
\newcommand{\Z}{\mathbb Z}
\newcommand{\N}{\mathbb N}
\newcommand{\K}{{\mathbb K}}
\renewcommand{\S}[2]{\textnormal{S}_{#1}\left(#2\right)}

\newcommand{\ie}{i.e.,\ }

\newcommand{\shuffle}{\, \raisebox{1.2ex}[0mm][0mm]{\rotatebox{270}{$\exists$}} \,}
\newcommand{\HarmonicSumsP}{\texttt{HarmonicSums}}
\def\f{\frac}
\def\bi{\binom}
\def\l{\left}
\def\r{\right}

\theoremstyle{plain}
\newtheorem{theorem}{Theorem}
\newtheorem{lemma}[theorem]{Lemma}

\theoremstyle{definition}

\newtheorem{example}[theorem]{Example}

\theoremstyle{remark}

\allowdisplaybreaks[4]

\title{\bf Discovering and Proving Infinite\\ Binomial Sums Identities}

\author{Jakob Ablinger\thanks{Supported by the Austrian Science Fund (FWF) grant SFB F50 (F5009-N15)}\\
\small Research Institute for Symbolic Computation\\[-0.8ex]
\small Johannes Kepler University\\[-0.8ex] 
\small Linz, Austria\\
}

\date{\dateline{July 29, 2015}}

\begin{document}

\maketitle

\begin{abstract}
We consider binomial and inverse binomial sums at infinity and rewrite them in terms of a small set of constants, such as powers of $\pi$ or $\log(2)$. In order to perform these simplifications, we view the series as specializations of generating series. 
For these generating series, we derive integral representations in terms of root-valued iterated integrals. Using substitutions, we express the iterated integrals as cyclotomic 
harmonic polylogarithms. Finally, by applying known relations among the cyclotomic harmonic polylogarithms, we derive expressions in terms of several constants.
\end{abstract}


\section{Introduction}
\label{sec:1}
The goal of this article is to find and prove identities of the following form: 
\begin{eqnarray*}
\sum_{i=1}^{\infty } \frac{\sum_{j=1}^i \frac{1}{j^2}}{(1+2 i) \binom{2 i}{i}}&=& \frac{\pi ^3}{81 \sqrt{3}}, \\
\sum_{i=1}^{\infty } \frac{3^i \sum_{j=1}^i \frac{1}{j^2}}{(1+2 i) \binom{2 i}{i}}&=& \frac{8 \pi ^3}{81 \sqrt{3}}, \\
\sum_{i=1}^{\infty } \frac{\sum_{j=1}^i \frac{1}{j^3}}{i^2 \binom{2 i}{i}}&=& \frac{\zeta_5}{9}+\frac{\pi ^2 \zeta_3}{27},\\
\sum_{i=1}^{\infty } \frac{1}{i^5 \binom{2 i}{i}}&=& \frac{9}{4} \sqrt{3} \pi  c_ 8-\frac{19 \zeta_5}{3}+\frac{\pi ^2 \zeta_3}{9}+\frac{9 \sqrt{3} \pi }{4}-\frac{\pi ^5}{27 \sqrt{3}},
\end{eqnarray*}
with $\zeta_n:=\sum_{i=1}^\infty\frac{1}{i^n}$ and $c_8:=\sum_{i=0}^\infty\frac{1}{(3\;i+1)^4}.$

Note that binomial and inverse binomial sums (see for example \cite{Ablinger:2014}) are of interest in physics: in particular, these sums have been studied in order to perform calculations of higher 
order corrections to scattering processes in particle physics \cite{Ablinger:2013eba,Ogreid:1997bx,Kalmykov:2000qe,Fleischer:1998nb,Davydychev:2001,Davydychev:2003mv,Jegerlehner:2003,Kalmykov:2007dk,Weinzierl:2004bn}. Central binomial sums were also considered in \cite{Lehmer:1985,Zucker:1985,Borwein:2001,Borwein:2000}, 
and there is a connection to Ap\'ery's proof of the irrationality of $\zeta(3)$ (see \cite{Borwein:1987}).
In \cite{ZhiWei:2014}, a list of conjectures on series for powers of $\pi$ and other important constants is presented. In the frame of this article we are going to prove several of these conjectures involving binomial and inverse binomial sums.
We summarize our approach with a concrete example.
Consider the sum 
\begin{equation}
\sum_{k=1}^\infty\frac{3\; \S{1}{k}-\frac{1}{k}}{k^2\binom{2k}{k}}, \label{examplesum}
\end{equation}
where $$\S{a}{k}:=\sum_{i=1}^n\frac{\sign{a}^k}{k^{\abs{a}}}$$
denotes the generalized harmonic numbers with $a\in\Z\setminus\{0\}.$
As a first step we derive an integral representation for (\ref{examplesum}). Using tools from \cite{Ablinger:2014} this leads to 
\begin{eqnarray}\label{examplesumintrep}
 \sum_{k=1}^\infty\frac{3\; \S{1}{k}-\frac{1}{k}}{k^2\binom{2k}{k}}&=&\frac{3}{2}\biggl( \int_0^1\frac{1}{t}\int_0^t\frac{1}{\sqrt{s(4-s)}}\int_0^s\frac{\sqrt{r}}{\sqrt{4-r}}drdsdt\nonumber\\
						    &&+\int_0^1\frac{1}{\sqrt{t(4-t)}}\int_0^t\frac{\sqrt{s}}{\sqrt{4-s}}dsdt\nonumber\\
						    &&+\int_0^1\frac{1}{\sqrt{t(4-t)}}\int_0^t\frac{1}{4-s}\int_0^s\frac{\sqrt{r}}{\sqrt{4-r}}drdsdt+1\biggr)\nonumber\\
						    &&-\int_0^1\frac{1}{t}\int_0^t\frac{1}{\sqrt{s(4-s)}}\int_0^s\frac{1}{\sqrt{r(4-r)}}drdsdt.
\end{eqnarray}
Next we want to rewrite these iterated integrals in terms of so called cyclotomic harmonic polylogarithms \cite{Ablinger:2013jta,Ablinger:2013hcp,Ablinger:2011te} which are extensions of the harmonic polylogarithms \cite{Remiddi:1999ew}. In order to define cyclotomic harmonic polylogarithms we introduce the following
auxiliary function: For $a \in \N$ and $b \in \N,$ $b < \varphi(a)$ (here $\varphi$ denotes Euler's totient function)  we define $f_a^b:(0,1)\mapsto \R$ by 
\begin{eqnarray}
&&f_a^b(x)=\left\{ 
		\begin{array}{ll}
				\frac{1}{x} &  \textnormal{if }a=b=0  \\
				\frac{x^b}{\Phi_a(x)} & \textnormal{otherwise}
		\end{array} 
		\right.  \nonumber
\end{eqnarray}
where $\Phi_a(x)$ denotes the $a$th cyclotomic polynomial.

Cyclotomic polylogarithms are now defined recursively: Let $m_i=(a_i,b_i) \in \N^2,$ $b_i<\varphi(a_i);$ we define for $x\in (0,1):$
\begin{eqnarray}\label{hpldef}
\H{}{x}&=&1,\nonumber\\
\H{m_1,\ldots,m_k}{x} &=&\left\{ 
		  	\begin{array}{ll}
						\frac{1}{k!}(\log{x})^k,&  \textnormal{if }m_i=(0,0) \textnormal{ for } 1\leq i \leq k\\
						  &\\
						\int_0^x{f_{a_1}^{b_1}(y) \H{m_2,\ldots,m_k}{y}dy},& \textnormal{otherwise}. 
			\end{array} \right.  
\end{eqnarray}
The length $k$ of the vector $\ve m=(m_1,\cdots,m_k)$ is called the \textit{weight} of the cyclotomic harmonic polylogarithm $\H{\ve m}x.$
If $m_1\neq (1,0)$ then the limit $\lim_{x\to 1} \H{m_1,\ldots,m_w}x$ is finite and we define $$\H{m_1,\ldots,m_w}1:=\lim_{x\to 1} \H{m_1,\ldots,m_w}x.$$

By using the substitution $x\to \frac{(u-1)^2}{1+u+u^2}$ we are going to show below how we can express the iterated integrals from equation~(\ref{examplesumintrep}) in terms of cyclotomic harmonic polylogarithms. Hence we get
\begin{eqnarray}
\sum_{k=1}^\infty\frac{3\; \S{1}{k}-\frac{1}{k}}{k^2\binom{2k}{k}}&=&\frac{1}{4} \bigg(
        18\, \H{(2,0)}{1} \H{(3,0)}{1}
        -9\, \H{(3,0)}{1}^2
        -18\, \H{(3,0)}{1} \H{(3,1)}{1}
        \nonumber\\&&-18\, \H{(2,0),(3,0)}{1}
        -48\, \H{(3,0)}{1} \H{(3,0),(1,0)}{1}
        -18\, \H{(3,0),(2,0)}{1}
        \nonumber\\&&+18\, \H{(3,0),(3,0)}{1}
        +48\, \H{(3,1)}{1} \H{(3,0),(3,0)}{1}
        +18\, \H{(3,0),(3,1)}{1}
        \nonumber\\&&+18\, \H{(3,1),(3,0)}{1}
        -48\, \H{(3,0)}{1} \H{(3,1),(3,0)}{1}
        +48\, \H{(3,0),(1,0),(3,0)}{1}
        \nonumber\\&&+72\, \H{(3,0),(2,0),(3,0)}{1}
        +48\, \H{(3,0),(3,0),(1,0)}{1}
        -12\, \H{(3,0),(3,0),(3,0)}{1}
        \nonumber\\&&-72\, \H{(3,0),(3,1),(3,0)}{1}
        +48\, \H{(3,1),(3,0),(3,0)}{1}
\bigg).\label{examplesumcyclorep}
\end{eqnarray}

Finally, by using known relations (see~\cite{Ablinger:2013hcp,Ablinger:2011te}) and new relations worked out in Section~4 between cyclotomic harmonic polylogarithms at~1 we can derive
\begin{eqnarray}
\sum_{k=1}^\infty\frac{3\; \S{1}{k}-\frac{1}{k}}{k^2\binom{2k}{k}}&=&\zeta_3.
\end{eqnarray}
Summarizing, the proposed strategy, which has been implemented in the Mathematica package {\tt HarmonicSums}\footnote{The package {\tt HarmonicSums} can be downloaded at\\ \url{http://www.risc.jku.at/research/combinat/software/HarmonicSums}.} \cite{HarmonicSums}, is
\begin{enumerate}
 \item Rewrite the sums in terms of nested integrals (see Section~2).
 \item Rewrite the integrals in terms of cyclotomic harmonic polylogarithms (see Section~3).
 \item Provide a sufficiently strong database to eliminate relations among these cyclotomic polylogarithms and find reduced expressions (see Section~4).
\end{enumerate}
Using this strategy the two main computational challanges are to transform the nested integrals into expressions in terms of cyclotomic harmonic polylogarithms, and to calculate the underlying relations such that the expressions in terms of cyclotomic harmonic 
polylogarithms collapse to simple constants. Here the computational effort grows exponentially with the nested depth of the integrals. 
We are able express a large set of binomial and inverse binomial sums at infinity in terms of several constants using the proposed method. As a bonus, we can generate base identities that can be combined to new interesting results. 
Note that the main purpose of this article is to present this method which can be automated, hence not all identities presented in this paper are new identities. Many of the identities can already be found in \cite{Lehmer:1985,Zucker:1985,Kalmykov:2000qe,Fleischer:1998nb,Davydychev:2001,Davydychev:2003mv,Jegerlehner:2003,Kalmykov:2007dk}.
However we are also able to find new identies and prove several conjectures. In particular we are able to discover and prove some of the conjectures from \cite{ZhiWei:2014}. 

The remainder of this article is organized as follows: In Section~2, we show in detail how we can derive integral representations for special binomial and inverse binomial sums. In Section~3, we show how these integral representations can be transformed to 
expressions in terms of cyclotomic harmonic polylogarithms. Section~4 deals with relations between the cyclotomic harmonic polylogarithms. Finally, in Section~5, we summarize the base identities that we 
found together with some nice combined results and list the conjectures from \cite{ZhiWei:2014} that we could prove using our approach.

\section{Generating Functions and Infinite Nested Binomial Sums}
\label{sec:2}

As a first step, we derive integral representations for the binomial sums. In order to accomplish this task, we view infinite sums as specializations of generating functions \cite{Ablinger:2014}. Namely, if we are given an
integral representation of the generating 
function of a sequence, then we can obtain an integral representation for the
infinite sum over that sequence if the limit $x \to 1$ can be carried out. This approach to infinite sums can be summarized by the
following formula:
\[
 \sum_{i=1}^\infty f(i)  = \lim_{x\to1}\sum_{i=1}^\infty x^if(i).
\]
In order to find integral representations of the generating functions, we rely on some results given in \cite{Ablinger:2014}.
First we need the following well known properties.
\begin{lemma} Let $\K$ be a field of characteristic 0 and let $f:\N \to \K.$ Then the following identities hold in the ring $\K[[x]]$ of formal power series
 \begin{eqnarray}
  \sum_{n=1}^\infty\frac{x^n}{n}f(n) &=& \int_0^x\frac{1}{t}\sum_{n=1}^\infty t^nf(n)dt,\label{eq:GenFunReciprocal}\\
  \sum_{n=1}^\infty x^n\sum_{i=1}^nf(i) &=& \frac{1}{1-x}\sum_{n=1}^\infty x^nf(n),\label{eq:GenFunSum}\\
  \sum_{n=1}^\infty\frac{x^n}{n+1}f(n) &=& \sum_{n=1}^\infty\frac{x^n}{n}f(n)-\frac{1}{x}\int_0^x\sum_{n=1}^\infty \frac{t^n}{n}f(n)dt.
 \end{eqnarray}
\end{lemma}

In addition, we make use of the following identities, which are useful for expressions involving binomial coefficients. 
Related formulae can also be found in the Appendix of \cite{Fleischer:1998nb}, which do not 
explicitly express the results as iterated integrals.

\begin{lemma} Let $\K$ be a field of characteristic 0 and let $f:\N \to \K.$ Then the following identities hold in the ring $\K[[x]]$ of formal power series:
 \begin{eqnarray}
   \sum_{n=1}^\infty x^n\binom{2n}{n}\sum_{i=1}^nf(i) &=& 
      \frac{1}{4\sqrt{\frac{1}{4}-x}}\int_0^x\frac{\sum_{n=1}^\infty t^nn\binom{2n}{n}f(n)}{t\sqrt{\frac{1}{4}-t}}dt,\label{eq:GenFunBinomialSum}\\
   \sum_{n=1}^\infty \frac{x^n}{n\binom{2n}{n}}\sum_{i=1}^nf(i) &=&
      \sum_{n=1}^\infty \frac{x^n}{n\binom{2n}{n}}f(n)+\sqrt{\frac{x}{4-x}}\int_0^x\frac{\sum_{n=1}^\infty \frac{t^n}{\binom{2n}{n}}f(n)}{\sqrt{t(4-t)}}dt, \label{eq:GenFunInverseBinomialSum1b}\\
   \sum_{n=1}^\infty\frac{x^n}{(2n+1)\binom{2n}{n}}\sum_{i=1}^nf(i) &=& 
      \frac{2}{\sqrt{x(4-x)}}\int_0^x\frac{\sum_{n=1}^\infty\frac{t^n}{\binom{2n}{n}}f(n)}{\sqrt{t(4-t)}}dt.\label{eq:GenFunInverseBinomialSum2}
 \end{eqnarray}
\end{lemma}


Let us illustrate the use of the formulae above by two simple examples.
\begin{example}
 Consider the generating function
 \begin{equation}\label{eq:GenFunExample1}
  \sum_{n=1}^\infty\frac{x^n}{n^3\binom{2n}{n}}.
 \end{equation}
 Applying \eqref{eq:GenFunReciprocal} twice we obtain
 \begin{eqnarray*}
  \sum_{n=1}^\infty\frac{x^n}{n^3\binom{2n}{n}} &=& \int_0^x\frac{1}{t}\sum_{n=1}^\infty \frac{t^n}{n^2\binom{2n}{n}}dt\\
    &=&\int_0^x\frac{1}{t}\int_0^t\frac{1}{s}\sum_{n=1}^\infty \frac{s^n}{n\binom{2n}{n}}dsdt.
 \end{eqnarray*}
 Now, by virtue of \eqref{eq:GenFunInverseBinomialSum1b} with $f(n):=\delta_{n,1},$ we obtain the result
 \begin{eqnarray}
  \sum_{n=1}^\infty\frac{x^n}{n^3\binom{2n}{n}} &=&\int_0^x\frac{1}{t}\int_0^t\frac{1}{\sqrt{s(4-s)}}\int_0^s\frac{1}{\sqrt{r(4-r)}}drdsdt.
 \end{eqnarray}
\end{example}

\begin{example}
 Consider the generating function
 \begin{equation}\label{eq:GenFunExample2}
  \sum_{n=1}^\infty\frac{x^n \S{1}n}{n^2\binom{2n}{n}}.
 \end{equation}
 Applying \eqref{eq:GenFunReciprocal} and then \eqref{eq:GenFunInverseBinomialSum1b}, we obtain
 \begin{eqnarray*}
  \sum_{n=1}^\infty\frac{x^n \S{1}n}{n^2\binom{2n}{n}} &=& \int_0^x\frac{1}{t}\sum_{n=1}^\infty\frac{t^n}{n^2\binom{2n}{n}}dt+\int_0^x\frac{\sqrt{t}}{\sqrt{4-t}}\int_0^t\frac{1}{\sqrt{s(4-s)}}\sum_{n=1}^\infty\frac{s^n}{n\binom{2n}{n}}dsdt.
 \end{eqnarray*}
 Again, by applying \eqref{eq:GenFunReciprocal} and then \eqref{eq:GenFunInverseBinomialSum1b} with $f(n):=\delta_{n,1},$ we obtain
 \begin{eqnarray}
  \sum_{n=1}^\infty\frac{x^n \S{1}n}{n^2\binom{2n}{n}} &=&\frac{1}{2}\biggl( \int_0^x\frac{1}{t}\int_0^t\frac{1}{\sqrt{s(4-s)}}\int_0^s\frac{\sqrt{r}}{\sqrt{4-r}}drdsdt\nonumber\\
						    &&+\int_0^x\frac{1}{\sqrt{t(4-t)}}\int_0^t\frac{\sqrt{s}}{\sqrt{4-s}}dsdt\nonumber\\
						    &&+\int_0^x\frac{1}{\sqrt{t(4-t)}}\int_0^t\frac{1}{4-s}\int_0^s\frac{\sqrt{r}}{\sqrt{4-r}}drdsdt+x\biggr).
 \end{eqnarray}
\end{example}

Note that combining Example \ref{eq:GenFunExample1} and \ref{eq:GenFunExample2} and setting $x\rightarrow 1$ we arrive at
 \begin{eqnarray}
  \lim_{x\to1}\left(3\;\sum_{n=1}^\infty\frac{x^n \S{1}n}{n^2\binom{2n}{n}}-\sum_{n=1}^\infty\frac{x^n}{n^3\binom{2n}{n}}\right)&=&\sum_{k=1}^\infty\frac{3\; \S{1}{k}-\frac{1}{k}}{k^2\binom{2k}{k}}.
 \end{eqnarray}
Since the limit $x\rightarrow 1$ can be computed, we obtain (\ref{examplesumintrep}) as our result.


\section{Special Iterated Integrals to Cyclotomic Harmonic Polylogarithms}
\label{sec:3}
In this section, we show how we can transform the iterated integrals that occur in the integral representation of Section \ref{sec:2} to 
expressions in terms of cyclotomic harmonic polylogarithms defined in (\ref{hpldef}). As an example, we consider the iterated integral 
\begin{equation}
 \int _0^1\int _0^t\int _0^s\frac{\sqrt{r}}{\sqrt{(4-t) t} (4-s) \sqrt{4-r}}drdsdt
\end{equation}
which is part of (\ref{examplesumintrep}). Substituting $r\to \frac{(u-1)^2}{1+u+u^2}$ in the innermost integral leads to
\begin{equation}
 \int _0^1\int _0^t\int _1^{\frac{-2-s+\sqrt{3} \sqrt{4 s-s^2}}{2 (-1+s)}}-\frac{\sqrt{3} (1-u)^2}{\sqrt{(4-t) t} (4-s)
   \left(1+u+u^2\right)^2}dudsdt.
\end{equation}
Proceeding by the substitution $s\to \frac{(v-1)^2}{1+v+v^2}$
\begin{equation}
 \int _0^1\int _1^{\frac{-2-t+\sqrt{3} \sqrt{4 t-t^2}}{2 (-1+t)}}\int _1^v-\frac{(1-v) \left(-\left(\sqrt{3}
   (1-u)^2\right)\right)}{\sqrt{(4-t) t} \left((1+v) \left(1+v+v^2\right)\right) \left(1+u+u^2\right)^2}dudvdt
\end{equation}
and finally by $t\to \frac{(w-1)^2}{1+w+w^2}$ we arrive at
\begin{equation}
 \int _1^0\int _1^w\int _1^v-\frac{(1-v)(1-u)^2}{\left(1+w+w^2\right) \left((1+v)
   \left(1+v+v^2\right)\right) \left(1+u+u^2\right)^2}dudvdw.
\end{equation}
Since the integration for cyclotomic harmonic polylogarithms always starts at 0, we rewrite this integral in the form
\begin{eqnarray}
 &&3 \int_0^1 \frac{1}{1+w+w^2} \, dw \int_0^1 \frac{1-v}{(1+v) \left(1+v+v^2\right)} \, dv \int_0^1\frac{(1-u)^2}{\left(1+u+u^2\right)^2} \, du\nonumber\\
 &&-3 \int_0^1 \frac{1}{1+w+w^2} \, dw \int _0^1\int _0^v\frac{(1-v)(1-u)^2}{\left((1+v) \left(1+v+v^2\right)\right) \left(1+u+u^2\right)^2}dudv\nonumber\\
 &&-3 \int _0^1\int _0^w\frac{1-v}{\left(1+w+w^2\right)\left((1+v) \left(1+v+v^2\right)\right)}dvdw \int_0^1 \frac{(1-u)^2}{\left(1+u+u^2\right)^2} \, du\nonumber\\
 &&+3 \int _0^1\int _0^w\int_0^v\frac{(1-v) (1-u)^2}{\left(1+w+w^2\right) \left((1+v) \left(1+v+v^2\right)\right) \left(1+u+u^2\right)^2}dudvdw
 \label{allmosthform}
\end{eqnarray}
These integrals can be rewritten in terms of cyclotomic harmonic polylogarithms using partial fractions and the formula
\begin{eqnarray}
\int_a^b\frac{f(x)}{g(x)^{i}}dx&=&\frac{1}{i-1}\left(\frac{f(b)}{g(b)^{i-1}g'(b)}-\frac{f(a)}{g(a)^{i-1}g'(a)}\right.\nonumber\\&&\left.+\int_a^b\frac{1}{g(x)^{i-1}}\left(\frac{f(x)g''(x)}{g'(x)^2}-\frac{f'(x)}{g'(x)}\right) dx\right)
\label{intdegred}
\end{eqnarray}
which can be derived by using integration by parts on 
$$
\int_a^b\frac{f'(x)}{g(x)^{i-1}g'(x)}dx.
$$
Let us illustrate this on 
\begin{eqnarray}
\int _0^1\int _0^v\frac{(1-v)(1-u)^2}{\left((1+v) \left(1+v+v^2\right)\right) \left(1+u+u^2\right)^2}dudv.
\label{intrewriteexample}
\end{eqnarray}
First we look at the inner integral
\begin{eqnarray}
 \int _0^v\frac{(1-u)^2}{\left(1+u+u^2\right)^2}du&=&\int _0^v\frac{1}{\left(1+u+u^2\right)^2}du-2\int _0^v\frac{u}{\left(1+u+u^2\right)^2}du\nonumber\\&&+\int _0^v\frac{u^2}{\left(1+u+u^2\right)^2}du
 \label{intrewriteexample1}
\end{eqnarray}
Applying (\ref{intdegred}) to the first integral twice, together with partial fractioning yields
\begin{eqnarray}
\int _0^v\frac{1}{\left(1+u+u^2\right)^2}du&=&\frac{1}{(2 v+1) \left(v^2+v+1\right)}-1+\int_0^v \frac{4}{3 (2 u+1)^2}-\frac{1}{3 \left(u^2+u+1\right)} du\nonumber\\
      &=&\frac{1}{(2 v+1) \left(v^2+v+1\right)}-1+\frac{4 v}{3 (2 v+1)}-\int_0^v\frac{1}{3 \left(u^2+u+1\right)} du\nonumber\\
      &=&-\frac{(-1+v) v}{3 \left(1+v+v^2\right)}+\frac{2}{3} \H{(3,0)}v.
\end{eqnarray}
Applying the second integral and using the same strategy yields
\begin{eqnarray*}
&&\int _0^1\frac{(1-v)}{\left(1+v \left(1+v+v^2\right)\right)}\left(-\frac{(-1+v) v}{3 \left(1+v+v^2\right)}+\frac{2}{3} \H{(3,0)}v\right)dv=-\frac{4}{3} \H{(2,0)}1+\H{(3,0)}1\\ 
&&\hspace{4cm}+\frac{4}{3} \H{(3,1)}1+\frac{4}{3} \H{(2,0),(3,0)}1-\frac{2}{3} \H{(3,0),(3,0)}1-\frac{4}{3} \H{(3,1),(3,0)}1.
\end{eqnarray*}
Proceeding in the same manner for the other two integrals in (\ref{intrewriteexample1}) and combining yields an expression in terms of cyclotomic harmonic polylogarithms for (\ref{intrewriteexample}):
\begin{eqnarray*}
&&4 \H{(2, 0), (3, 0)}1 - 2 \H{(3, 0), (3, 0)}1 - 4 \H{(3, 1), (3, 0)}1\\ &&- 2 \H{(2, 0)}1 + 3 \H{(3, 0)}1 + 2 \H{(3, 1)}1-1.
\end{eqnarray*}

In a similar fashion we can rewrite all the iterated integrals in (\ref{allmosthform}) in terms of cyclotomic harmonic polylogarithms. This leads to
\begin{eqnarray}
&&12 \H{(3,0),(2,0),(3,0)}{1}-6 \H{(3,0),(3,0),(3,0)}{1}
-12 \H{(3,0),(3,1),(3,0)}{1}+\big(
        12 \H{(3,0),(3,0)}{1}
        \nonumber\\&&+12 \H{(3,0),(3,1)}{1}
        +12 \H{(3,1),(3,0)}{1}
        -12 \H{(2,0),(3,0)}{1}
         -12 \H{(3,0),(2,0)}{1}     
\big) \H{(3,0)}{1}\nonumber\\
&&+\big(
        12 \H{(2,0)}{1}
        -12 \H{(3,1)}{1}
        -6
\big) \H{(3,0)}{1}^2
-6 \H{(3,0)}{1}^3
+6 \H{(3,0),(3,0)}{1}+1.
\end{eqnarray}
Note that the computational effort that has to be made grows exponentially with the nested depth of the integrals. In \HarmonicSumsP\ sophisticated algorithms are used to speed up these computations.
By applying these substitutions to all iterated integrals in~(\ref{examplesumintrep}) we can derive~(\ref{examplesumcyclorep}).

Note that in our computations, integrands of the form $\sqrt{4-3 x}\sqrt{x}$ and $\sqrt{2-x}\sqrt{x}$ also appear. In order to deal with these integrands, 
we used the substitutions $x\to\frac{-1 + u^2}{-1 - \sqrt{2}u + u^2}$ and $x\to\frac{(1-u)^2}{1+u^2},$ respectively. More precisely we utilize the following substitutions:

 \begin{eqnarray*}
   &\sqrt{(4-x)x}:& x\to\frac{(1-u)^2}{1+u+u^2} \\ 
   &\sqrt{(4-3 x)x}:& x\to\frac{-1 + u^2}{-1 - \sqrt{2}u + u^2} \\ 
   &\sqrt{(2-x)x}:& x\to\frac{(1-u)^2}{1+u^2} \\ 
   &\sqrt{(1-x)x}:& x\to\frac{(1-u)^2}{1+u^2}\\ 
   &\sqrt{1-x}:& x\to\frac{-4u}{(1-u)^2}.
 \end{eqnarray*}


\section{Relations between Cyclotomic Harmonic Polylogarithms}
\label{sec:4}
In this section, we examine relations beween cyclotomic harmonic polylogarithms at~1 (compare \cite{Ablinger:2013hcp,Ablinger:2011te,Ablinger:2013cf}). A first set of relations originates from the shuffle algebra structure of 
cyclotomic harmonic polylogarithms \ie the product of two cyclotomic harmonic polylogarithms of the same argument can be expressed using the formula
\begin{equation}
\label{CShpro}
\H{\ve p}x\H{\ve q}x=\sum_{\ve r= \ve p \shuffle \ve q}\H{\ve r}x
\end{equation}
in which $\ve p \shuffle \ve q$ represents all merges of $\ve p$ and $\ve q$ in which the relative orders of the elements of $\ve p$ and $\ve q$ are preserved. As an example, we have:
 \begin{eqnarray*}
\text{H}_{(3,0),(2,0)}(1) \text{H}_{(3,1),(3,0)}(1)&=&\text{H}_{(3,0),(2,0),(3,1),(3,0)}(1)+\text{H}_{(3,0),(3,1),(2,0),(3,0)}(1)\\&&+\text{H}_{(3,0),(3,1),(3,0),(2,0)}(1)+\text{H}_{(3,1),(3,0),(2,0),(3,0)}(1)\\&&+2 \text{H}_{(3,1),(3,0),(3,0),(2,0)}(1).
 \end{eqnarray*}
Secondly, there are the so-called duality relations which are due to argument transforms of the form $1-x\to x$ or $x\to \frac{a-x}{b+x},\; a,b\in\R^*.$ The transform $x\to \frac{1-x}{1+x}$ leads, for example, to the relation
$$\text{H}_{(0,0),(1,0)}(1)-\text{H}_{(0,0),(2,0)}(1)+\text{H}_{(1,0),(0,0)}(1)-\text{H}_{(2,0),(1,0)}(1)+\text{H}_{(2,0),(2,0)}(1)=0.$$
Other classes of relations for cyclotomic harmonic polylogarithm at~1 originate from the sum representation of these constants \cite{Ablinger:2011te}.
For $a_i,k \in \N^*,$ $b_i,n \in \N$ and $c_i\in\Z^*$ we define
\begin{eqnarray*}
\S{(a_1,b_1,c_1),\ldots,(a_k,b_k,c_k)}{n}=\sum_{n \geq i_1 \geq \cdots \geq i_k \geq 1}\frac{\sign{c_1}^{i_1}}{(a_1 i_1+b_1)^{\abs{c_1}}}\cdots\frac{\sign{c_k}^{i_k}}{(a_k i_k+b_k)^{\abs{c_k}}};
\end{eqnarray*}
$k$ is called the \textit{depth} and $w=\sum_{i=1}^k\abs{c_i}$ is called the \textit{weight} of the cyclotomic harmonic sum $\S{(a_1,b_1,c_1),\ldots,(a_k,b_k,c_k)}{n}$.

Further relations can be discovered by utilizing the connection between cyclotomic harmonic polylogarithms at~1 and cyclotomic harmonic sums at~$\infty.$
This link can be established via the power series expansion of cyclotomic harmonic polylogarithms. In the following we sketch this translation mechanism, for details 
we refer the reader to \cite{Ablinger:2013hcp,Ablinger:2011te}.
Let $\Phi_a(x)$ be a cyclotomic polynomial. Then we can write 
$${\frac{1}{\Phi_a(x)}}=\sum_{q=0}^{a-1}f_q\sum_{i=0}^{\infty}x^{a i +q}.$$
Hence we get the power series expansion of depth one cyclotomic harmonic polylogarithms as follows.
\begin{eqnarray*}
\H{(a,b)}x&=&\int_0^x{\frac{y^b}{\Phi_a(y)}}dy=\sum_{q=0}^{a-1}f_q\int_0^x\sum_{i=0}^{\infty}y^{a i+q+b}\\
	    &=&\sum_{q=0}^{a-1}f_q\sum_{i=0}^{\infty}\frac{ x^{a i+q+b+1}}{a i+q+b+1}\\
	    &=&\sum_{q=0}^{a-1}f_q\sum_{i=1}^{\infty}\frac{x^{a i+q-a+b+1}}{a i+q-a+b+1}.
\end{eqnarray*}
We can proceed recursively on the depth to obtain a power series expansion of a general cyclotomic harmonic polylogarithm. 
Let $\H{\ve m}x$ be a cyclotomic harmonic polylogarithm with
$$\H{\ve m}x=\sum_{j=1}^w\sum_{i=1}^{\infty}\frac{x^{z i+c_j}}{(z i+c_j)^{g_j}}\S{\ve n_j}i$$
for $x \in (0,1)$, $w,g_j\in \N, c_j\in \Z$  and some cyclotomic harmonic sums $\S{\ve n_j}i.$\\ 
Provided that $a=0$ or $a k=z$ for some $k \in \N$ we have
\begin{eqnarray*}
\H{(0,0),\ve m}x&=&\sum_{j=1}^w\sum_{i=1}^{\infty}\frac{x^{z i+c_j}}{(z i+c_j)^{g_j+1}}\S{\ve n_j}i,\\
\H{(a,b),\ve m}x&=&\sum_{j=1}^w{\sum_{q=0}^{z-1}f_q\sum_{i=1}^{\infty}\frac{ x^{z i+q+b+c_j+1}}{(z i+q+b+c_j+1)}\S{(z,c_j,g_j),\ve n_j}i}
\end{eqnarray*}
where $$\frac{1}{\Phi_a(x)}=\sum_{q=0}^{z-1}f_q\sum_{p=0}^{\infty}x^{z p +q}.$$
Finally, for $x\rightarrow1$ these sums turn into cyclotomic harmonic sums at infinity if $c\neq1$:
$$
\sum_{i=1}^\infty x^{ai+b} \frac{ \S{\ve n}i}{(ai+b)^c}\rightarrow \S{(a,b,c),\ve n}{\infty}.
$$
\begin{example} Consider
\begin{eqnarray*}
\textnormal{H}_{(0,0),(3,2)}(x)&=& \frac{1}{9} \bigg(
        \sum_{\tau_1=1}^{\infty } \frac{\big(
                x^3\big)^{\tau_1}}{\tau_1^2}
        -9 x 
        \sum_{\tau_1=1}^{\infty } \frac{\big(
                x^3\big)^{\tau_1}}{\big(
                1+3 \tau_1\big)^2}
\bigg).
\end{eqnarray*}
Then, for $x\rightarrow1$ we get 
$$
\H{(0,0),(3,2)}{1}=\frac{1}{9}\S{(1,0,2)}{\infty}-\S{(3,1,2)}{\infty}.
$$
\end{example}

In summary we can rewrite a cyclotomic harmonic polylogarithm at~1 in terms of cyclotomic harmonic sums at~$\infty$ using the power series expansion. From the sum representations, we can deduce additional classes of relations:
cyclotomic harmonic sums form a quasi-shuffle algebra (see \cite{Ablinger:2013hcp,Ablinger:2011te,Hoffman}) and hence we get relations of the form
\begin{eqnarray*}
  &&\;\textnormal{S}_{(5,3,3)}(\infty) \;\textnormal{S}_{(2,1,2),(3,2,1)}(\infty)
  =120 \;\textnormal{S}_{(2,1,1),(3,2,1)}(\infty)+19 \;\textnormal{S}_{(2,1,2),(3,2,1)}(\infty)\\&&\hspace{2cm}
  -45 \;\textnormal{S}_{(2,1,2),(5,3,1)}(\infty)+15 \;\textnormal{S}_{(2,1,2),(5,3,2)}(\infty)
  -5 \;\textnormal{S}_{(2,1,2),(5,3,3)}(\infty)\\&&\hspace{2cm}-300 \;\textnormal{S}_{(5,3,1),(3,2,1)}(\infty)
  -100 \;\textnormal{S}_{(5,3,2),(3,2,1)}(\infty)-25\;\textnormal{S}_{(5,3,3),(3,2,1)}(\infty)\\&&\hspace{2cm}
  +\textnormal{S}_{(2,1,2),(3,2,1),(5,3,3)}(\infty)+\textnormal{S}_{(2,1,2),(5,3,3),(3,2,1)}(\infty)
  +\textnormal{S}_{(5,3,3),(2,1,2),(3,2,1)}(\infty).                                                  	
\end{eqnarray*}
In addition to the quasi-shuffle relation there are two duplications: For $c_i \geq 1$ we have
\begin{eqnarray*}
\sum{\S{(a_m,b_m,\pm c_m),\ldots,(a_1,b_1, \pm c_1)}{2n}}=2^m \S{(2 a_m,b_m, c_m),\ldots,(2 a_1,b_1, c_1)}{n}
\end{eqnarray*}
where we sum on the left hand side over the $2^m$ possible combinations.\\
In addition, let $d_i \in \{-1,1\}.$ Then we have
\begin{eqnarray*}
\sum{d_m \cdots d_1 \S{(a_m,b_m,d_m c_m),\ldots,(a_1,b_1,d_1 c_1)}{2n}}=2^m \S{(2 a_m,b_m-a_m, c_m),\ldots,(2 a_1,b_1-a_1, c_1)}{n},
\end{eqnarray*}
where we sum on the left hand side over the $2^m$ possible combinations of the $d_i$.\\
Finally, there is the following multiple argument relation: 
For $a, k \in \N, b\in\N_0$, $c\in \Z^*$, $k\geq 2,$
\begin{eqnarray*}
\S{(a,b,c)}{k\; n}=\sum_{i=0}^{k-1}\sign{c}^i\S{(k\; a,b-a \, i, \sign{c}^k\abs{c})}{n};
\end{eqnarray*}
and for $a_i, m \in \N, b_i, k \in \N_0$, $c_i\in \Z^*$, $k\geq 2,$
\begin{eqnarray*}
&&\S{(a_m,b_m,c_m),(a_{m-1},b_{m-1}, c_{m-1}),\ldots,(a_1,b_1,c_1)}{k \; n}=\\
&&\hspace{1cm}\sum_{i=0}^{m-1}\sum_{j=1}^{n} \frac{\S{(a_{m-1},b_{m-1},c_{m-1}),\ldots,(a_1,b_1,c_1)}{k \; j-i}\sign{c_m}^{k\; j -i}} {(a_m (k\; j-i)+b_1)^{\abs{c_m}}}.
\end{eqnarray*}
As an example of a multiple argument relation, we state
\begin{eqnarray*}
\S{(2,1,2)}{3 n}&=&\frac{1}{9} \S{(2,1,2)}{n}+\S{(6,1,2)}{n}+\S{(6,5,2)}{n}-\frac{1}{(6 n+3)^2}-\frac{1}{(6 n+5)^2}+\frac{34}{225}.
\end{eqnarray*}
Note that if the $\lim_{n\to\infty}\S{\ve m}{n}$ exists then it makes no difference whether we consider sums at $\infty$ or $k\cdot\infty,\ k\in\N,$ and hence these relations remain valid in the limit. For example we get
\begin{eqnarray*}
\S{(2,1,2)}{\infty}&=&\frac{1}{9} \S{(2,1,2)}{\infty}+\S{(6,1,2)}{\infty}+\S{(6,5,2)}{\infty}+\frac{34}{225}.
\end{eqnarray*}
Summarizing we get the following classes of relations \cite{Ablinger:2011te}:
\begin{itemize}
\item stuffle relations (quasi shuffle algebra of cyclotomic sums);
\item two duplication relations (if the sum is finite, it makes no difference whether the argument is $\infty$ or $2\cdot\infty$);
\item multiple argument relations (if the sum is finite, it makes no difference whether the argument is $\infty$ or $k\cdot\infty$);
\item shuffle relations (shuffle algebra of cyclotomic polylogarithms);
\item duality relations of cyclotomic polylogarithms.
\end{itemize}

As a first example, we consider cyclotomy 4, \ie the cyclotomic polylogarithms with letters in $$\{(0, 0), (1, 0), (2, 0), (4, 0), (4, 1)\}$$ 
or the cyclotomic sums with letters in $$\left\{\frac{(\pm1)^i}{i},\frac{(\pm1)^i}{2i+1}\right\}$$ up to weight 5. Using the relations from above we can express all these constants using the following basis constants:
\begin{eqnarray*}
&&\H{(1,0)}{1},\H{(2,0)}{1},\H{(4,0)}{1};\\
&&\H{(0,0),(4,1)}{1};\\
&&\H{(1,0),(0,0),(0,0)}{1},\H{(0,0),(4,1),(4,0)}{1};\\
&&\H{(2,0),(1,0),(1,0),(1,0)}{1},\H{(4,0),(0,0),(0,0),(0,0)}{1},\H{(0,0),(4,1),(4,1),(4,0)}{1};\\
&&\H{(1,0),(0,0),(0,0),(0,0),(0,0)}{1},\H{(2,0),(1,0),(1,0),(1,0),(1,0)}{1},\H{(4, 1), (1, 0), (1, 0), (4, 0), (1, 0)}{1},\\
&&\H{(4, 1), (1, 0), (4, 0), (1, 0), (1, 0)}{1},\H{(4, 1), (4, 0), (1, 0), (1, 0), (1, 0)}{1},\H{(4, 1), (4, 1), (1, 0), (1, 0), (1, 0)}{1}.
\end{eqnarray*}

\begin{table}[t]\centering
\begin{tabular}{|c|r|r|r|r|r|}
\hline
\multicolumn{1}{|c}{\sf w}             &
\multicolumn{1}{|c}{1}                 &
\multicolumn{1}{|c}{2}                 &
\multicolumn{1}{|c}{3}                 &
\multicolumn{1}{|c|}{4}                &
\multicolumn{1}{|c|}{5}                \\
\hline
{\sf \# sums} 
&  4  & 20 & 100 & 500 & 2500\\
{\sf \# logs} 
&  5  & 25 & 125 & 625 & 3125\\
{\sf \# basis constants} 
&  3  &  1 &   2 &   3 &  6\\
\hline
\end{tabular}
\caption{\label{CShreltab4}Number of cyclotomic harmonic polylogarithms and cyclotomic harmonic sums together with the number of basis constants of cyclotomy 4 at different weights.}
\end{table}

As a second example we consider cyclotomy 6, \ie the cyclotomic polylogarithms with letters in $$\{(0, 0), (1, 0), (2, 0), (3, 0), (3, 1), (6, 0), (6, 1)\}$$ 
or the cyclotomic sums with letters in $$\left\{\frac{(\pm1)^i}{i},\frac{(\pm1)^i}{3i+1},\frac{(\pm1)^i}{3i+2}\right\}$$ up to weight 5. Using the relations from above we can express all these constants using the following 
basis constants (note that we only list the basis constants up to weight 4):
\begin{eqnarray*}
&&\H{(1,0)}{1},\H{(2,0)}{1},\H{(3,1)}{1},\H{(6,1)}{1};\\
&&\H{(0,0),(6,1)}{1},\H{(6,1),(3,1)}{1};\\
&&\H{(0,0),(0,0),(1,0)}{1},\H{(0,0),(3,1),(6,1)}{1},\H{(0,0),(6,1),(3,1)}{1},\H{(0,0),(6,1),(6,0)}{1},\H{(0,0),(6,1),(6,1)}{1},\\&&\H{(6,1),(6,1),(2,0)}{1},\H{(6,1),(6,1),(3,0)}{1};\\
&&\H{(0,0),(0,0),(0,0),(6,1)}{1},\H{(0,0),(0,0),(6,1),(3,1)}{1},\H{(0,0),(0,0),(6,1),(6,0)}{1},\H{(0,0),(0,0),(6,1),(6,1)}{1},\\&&
\H{(0,0),(6,1),(3,1),(3,1)}{1},\H{(0,0),(6,1),(3,1),(6,0)}{1},\H{(0,0),(6,1),(3,1),(6,1)}{1},\H{(0,0),(6,1),(6,0),(3,1)}{1},\\&&
\H{(0,0),(6,1),(6,0),(6,1)}{1},\H{(0,0),(6,1),(6,1),(2,0)}{1},\H{(0,0),(6,1),(6,1),(3,0)}{1},\H{(0,0),(6,1),(6,1),(3,1)}{1},\\&&
\H{(0,0),(6,1),(6,1),(6,0)}{1},\H{(0,0),(6,1),(6,1),(6,1)}{1},\H{(2,0),(1,0),(1,0),(1,0)}{1},\H{(6,1),(6,1),(3,1),(3,1)}{1},\\&&
\H{(6,1),(6,1),(6,1),(2,0)}{1},\H{(6,1),(6,1),(6,1),(3,1)}{1}.
\end{eqnarray*}

In the following we give some example relations for the cyclotomic harmonic polylogarithms appearing in (\ref{examplesumcyclorep}):
\begin{eqnarray*}
 \H{(3,0)}{1}&=&\H{(6,1)}{1}, \\
 \H{(2,0),(3,0)}{1}&=&\H{(2,0)}{1} \H{(6,1)}{1}
-
\frac{1}{4} \H{(6,1)}{1}^2
-\frac{1}{3} \H{(0,0),(6,1)}{1},
 \\
 \H{(3,0),(1,0)}{1}&=&\H{(3,1)}{1} \H{(6,1)}{1}
-\frac{1}{2} \H{(6,1)}{1}^2
-\frac{4}{3} \H{(0,0),(6,1)}{1},
 \\
 \H{(3,0),(2,0)}{1}&=&\frac{1}{4} \H{(6,1)}{1}^2
+\frac{1}{3} \H{(0,0),(6,1)}{1},
 \\
 \H{(3,0),(3,0)}{1}&=&\frac{1}{2} \H{(6,1)}{1}^2, \\
 \H{(3,0),(3,1)}{1}&=&\H{(3,1)}{1} \H{(6,1)}{1}
-\frac{1}{3} \H{(0,0),(6,1)}{1},
 \\
 \H{(3,1),(3,0)}{1}&=&\frac{1}{3} \H{(0,0),(6,1)}{1}, \\
 \H{(3,0),(1,0),(3,0)}{1}&=&-\frac{1}{2} \H{(6,1)}{1}^3
-\frac{2}{3} \H{(6,1)}{1} \H{(0,0),(6,1)}{1}
-\frac{2}{27} \H{(0,0),(0,0),(1,0)}{1},
 \\
 \H{(3,0),(2,0),(3,0)}{1}&=&-\frac{7}{4} \H{(6,1)}{1}^3
-\frac{7}{3} \H{(6,1)}{1} \H{(0,0),(6,1)}{1}
-\frac{13}{18} \H{(0,0),(0,0),(1,0)}{1},
 \\
 \H{(3,0),(3,0),(1,0)}{1}&=&\frac{1}{2} \H{(3,1)}{1} \H{(6,1)}{1}^2
-\frac{1}{3} \H{(6,1)}{1} \H{(0,0),(6,1)}{1}
+\frac{1}{27} \H{(0,0),(0,0),(1,0)}{1},
 \\
 \H{(3,0),(3,0),(3,0)}{1}&=&\frac{1}{6} \H{(6,1)}{1}^3, \\
 \H{(3,0),(3,1),(3,0)}{1}&=&-\frac{4}{3} \H{(6,1)}{1}^3
-\frac{5}{3} \H{(6,1)}
{1} \H{(0,0),(6,1)}{1}
-\frac{14}{27} \H{(0,0),(0,0),(1,0)}{1},
 \\
 \H{(3,1),(3,0),(3,0)}{1}&=&\frac{2}{3} \H{(6,1)}{1}^3
+\H{(6,1)}{1} \H{(0,0),(6,1)}{1}
+\frac{7}{27} \H{(0,0),(0,0),(1,0)}{1}.
 \\
\end{eqnarray*}

\begin{table}[t]\centering
\begin{tabular}{|c|r|r|r|r|r|}
\hline
\multicolumn{1}{|c}{\sf w}             &
\multicolumn{1}{|c}{1}                 &
\multicolumn{1}{|c}{2}                 &
\multicolumn{1}{|c}{3}                 &
\multicolumn{1}{|c|}{4}                &
\multicolumn{1}{|c|}{5}                \\
\hline
{\sf \# sums} 
&  6  & 42 & 294 & 2058& 14406\\
{\sf \# logs} 
&  7  & 49 & 343 & 2401& 16807\\
{\sf \# basis constants} 
&  4  &  2 &   7 &   18 &  52\\
\hline
\end{tabular}
\caption{\label{CShreltab}Number of cyclotomic harmonic polylogarithms and cyclotomic harmonic sums together with the number of basis constants of cyclotomy 6 at different weights.}
\end{table}

Plugging these relations into (\ref{examplesumcyclorep}) we find
\begin{eqnarray*}
\sum_{k=1}^\infty\frac{3\; \S{1}{k}-\frac{1}{k}}{k^2\binom{2k}{k}}=-\H{(0, 0), (0, 0), (1, 0)}1=\S{(1,0,3)}{\infty}=\zeta_3.
\end{eqnarray*}


\section{Results}
\label{sec:5}
In this section, we provide a number of base identities we discovered using the techniques outlined in Sections 2, 3, and 4. These results were obtained and proved using our implementation in the Mathematica package {\tt HarmonicSums} \cite{HarmonicSums}. Note that 
this section should illustrate what can be achieved using our implementation. Many of these base identities can already be found in \cite{Lehmer:1985,Zucker:1985,Kalmykov:2000qe,Fleischer:1998nb,Davydychev:2001,Davydychev:2003mv,Jegerlehner:2003,Kalmykov:2007dk}.
We define here a number of constants that appear. Note that these constants do not possess any further relations induced by the algebraic properties given in Section 4. 

\begin{tabular}{lll}
$C:= \text{Catalan};$ & $c_ 1:= \sum_{i=1}^{\infty } \frac{1}{(1+3 i)^2};$ & $c_ 2:= \sum_{i=1}^{\infty } \frac{(-1)^i \sum_{j=1}^i \frac{(-1)^j}{(1+3 j)^2}}{i};$\\
$c_ 3:= \sum_{i=1}^{\infty } \frac{\sum_{j=1}^i \frac{(-1)^j}{1+3 j}}{i^2};$&$c_ 4:= \sum_{i=1}^{\infty } \frac{\sum_{j=1}^i \frac{1}{1+3 j}}{i^2};$ & $c_ 5:= \sum_{i=1}^{\infty } \frac{(-1)^i \sum_{j=1}^i \frac{(-1)^j}{j}}{(1+2 i)^2};$ \\
$c_ 6:= \sum_{i=1}^{\infty } \frac{\sum_{j=1}^i \frac{(-1)^j}{2+3 j}}{(1+3 i)^2};$ & $c_ 7:= \sum_{i=1}^{\infty } \frac{(-1)^i}{(1+2 i)^4};$ & $c_ 8:= \sum_{i=1}^{\infty } \frac{1}{(1+3 i)^4};$ \\
$c_ 9:= \sum_{i=1}^{\infty } \frac{\sum_{j=1}^i \frac{(-1)^j}{(1+3 j)^2}}{i^2};$ & $c_ {10}:= \sum_{i=1}^{\infty } \frac{\sum_{j=1}^i \frac{(-1)^j}{1+3 j}}{i^3};$ & $c_ {11}:= \sum_{i=1}^{\infty } \frac{\sum_{j=1}^i \frac{1}{1+3 j}}{i^3};$ \\
$c_ {12}:=\sum_{i=1}^{\infty } \frac{1}{(1+3 i)^6};$ & $l_ 1:= \log (2);$ & $l_ 2:= \log (3);$\\
$l_ 3:= \log \big(7-4 \sqrt{3}\big);$ & $l_ 4:= \log \big(2-\sqrt{3}\big);$ & $l_ 5:= \log \big(2+\sqrt{3}\big);$\\
$p_ 1:= \text{Li}_ 2\left(\frac{1}{4}\right);$ & $p_ 2:= \text{Li}_ 2\left(1-\frac{\sqrt{3}}{2}\right);$ & $p_ 3:= \text{Li}_ 2\left(\frac{1}{4} \big( 2+\sqrt{3}\big)\right);$\\
$p_ 4:= \text{Li}_ 4\left(\frac{1}{2}\right);$ & \\
\\
\end{tabular}

Weight 1:
\begin{eqnarray}
\sum_{i=1}^{\infty } \frac{1}{i \binom{2 i}{i}}&=& \frac{\pi }{3 \sqrt{3}} \label{w5rel1}\\
\sum_{i=1}^{\infty } \frac{2^i}{i \binom{2 i}{i}}&=& \frac{\pi }{2}\\
\sum_{i=1}^{\infty } \frac{3^i}{i \binom{2 i}{i}}&=& \frac{2 \pi}{\sqrt{3}} \\
\sum_{i=1}^{\infty } \frac{1}{(1+2 i) \binom{2 i}{i}}&=& \frac{2 \pi}{3 \sqrt{3}}-1\\
\sum_{i=1}^{\infty } \frac{2^i}{(1+2 i) \binom{2 i}{i}}&=& \frac{\pi}{2}-1 \\
\sum_{i=1}^{\infty } \frac{3^i}{(1+2 i) \binom{2 i}{i}}&=& \frac{4 \pi }{3 \sqrt{3}}-1 \\
\sum_{i=1}^{\infty } \frac{16^{-i} \binom{2 i}{i}}{i}&=& 4 l_ 1-2 l_ 5\\
\sum_{i=1}^{\infty } \frac{16^{-i} \binom{2 i}{i}}{1+2 i}&=& \frac{\pi }{3}-1
\end{eqnarray}
Weight 2:
\begin{eqnarray}
\sum_{i=1}^{\infty } \frac{1}{i^2 \binom{2 i}{i}}&=& \frac{\pi ^2}{18} \label{w5rel2}\\
\sum_{i=1}^{\infty } \frac{
\sum_{j=1}^i \frac{1}{j}}{i \binom{2 i}{i}}&=& 2 c_ 1
-\frac{\pi  l_ 2}{3 \sqrt{3}}
+2
-\frac{5 \pi ^2}{54}
 \\
\sum_{i=1}^{\infty } \frac{2^i}{i^2 \binom{2 i}{i}}&=& \frac{\pi ^2}{8} \\
\sum_{i=1}^{\infty } \frac{2^i 
\sum_{j=1}^i \frac{1}{j}}{i \binom{2 i}{i}}&=& 2 C-\frac{\pi  l_ 1}{2}+\frac{\pi ^2}{8} \\
\sum_{i=1}^{\infty } \frac{3^i}{i^2 \binom{2 i}{i}}&=& \frac{2 \pi ^2}{9} \\
\sum_{i=1}^{\infty } \frac{3^i 
\sum_{j=1}^i \frac{1}{j}}{i \binom{2 i}{i}}&=& 9 c_ 1-\frac{4 \pi ^2}{9}+9 \\
\sum_{i=1}^{\infty } \frac{1}{(1+2 i)^2 \binom{2 i}{i}}&=& \frac{8 C}{3}-\frac{\pi  l_ 5}{3}-1 \\
\sum_{i=1}^{\infty } \frac{
\sum_{j=1}^i \frac{1}{j}}{(1+2 i) \binom{2 i}{i}}&=& 4 c_ 1
-\frac{2 \pi  l_ 2}{3 \sqrt{3}}
+4
-\frac{8 \pi ^2}{27}
 \\
\sum_{i=1}^{\infty } \frac{2^i 
\sum_{j=1}^i \frac{1}{j}}{(1+2 i) \binom{2 i}{i}}&=& 2 \
C-\frac{\pi  l_ 1}{2} \\
\sum_{i=1}^{\infty } \frac{3^i}{(1+2 i)^2 \binom{2 i}{i}}&=& 5 c_ 1
-\frac{\pi  l_ 2}
{3 \sqrt{3}}
+4
-
\frac{10 \pi ^2}{27}
 \\
\sum_{i=1}^{\infty } \frac{3^i 
\sum_{j=1}^i \frac{1}{j}}{(1+2 i) \binom{2 i}{i}}&=& 6 c_ 1-\frac{4 \pi ^2}{9}+6 \\
\sum_{i=1}^{\infty } \frac{
\sum_{j=1}^i \frac{1}{1+2 j}}{i \binom{2 i}{i}}&=& \frac{5 c_ 1}{2}
-\frac{\pi  l_ 2}{6 \sqrt{3}}
+\frac{9}{2}
-\frac{4 \pi }{3 \sqrt{3}}
-\frac{5 \pi ^2}{27}
 \\
\sum_{i=1}^{\infty } \frac{2^i}{i^2 \binom{2 i}{i}}&=& \frac{\pi ^2}{8} \\
\sum_{i=1}^{\infty } \frac{2^i 
\sum_{j=1}^i \frac{1}{1+2 j}}{i \binom{2 i}{i}}&=& 2 C+2-\pi  \\
\sum_{i=1}^{\infty } \frac{3^i}{i^2 \binom{2 i}{i}}&=& \frac{2 \pi ^2}{9} \\
\sum_{i=1}^{\infty } \frac{3^i 
\sum_{j=1}^i \frac{1}{1+2 j}}{i \binom{2 i}{i}}&=& \frac{15 c_ 1}{2}
+\frac{\pi  l_ 2}{\sqrt{3}}
+\frac{19}{2}
-\frac{8 \pi }{3 \sqrt{3}}
-\frac{5 \pi ^2}{9}
 \\
\sum_{i=1}^{\infty } \frac{
\sum_{j=1}^i \frac{1}{1+2 j}}{(1+2 i) \binom{2 i}{i}}&=& 5 c_ 1
-\frac{\pi  l_ 2}{3 \sqrt{3}}
+5
-\frac{2 \pi }{3 \sqrt{3}}
-\frac{10 \pi ^2}{27}
 \\
\sum_{i=1}^{\infty } \frac{2^i 
\sum_{j=1}^i \frac{1}{1+2 j}}{(1+2 i) \binom{2 i}
{i}}&=& 2 C-
\frac{\pi }{2} \\
\sum_{i=1}^{\infty } \frac{3^i 
\sum_{j=1}^i \frac{1}{1+2 j}}{(1+2 i) \binom{2 i}{i}}&=& 5 c_ 1
+\frac{2 \pi  l_ 2}{3 \sqrt{3}}
+5
-\frac{4 \pi }{3 \sqrt{3}}
-\frac{10 \pi ^2}{27}
 \\
\sum_{i=1}^{\infty } \frac{16^{-i} \binom{2 i}{i}}{i^2}&=& -12 l_ 1^2
+4 l_ 1 l_ 4
+8 l_ 1 l_ 5
-2 l_ 4 l_ 5
-l_5^2
-2 p_ 3
+\frac{\pi ^2}{3}
 \\
\sum_{i=1}^{\infty } \frac{16^{-i} \binom{2 i}{i} 
\sum_{j=1}^i \frac{1}{j}}{i}&=& 4 l_ 1^2
+l_2 l_ 3
+4 l_ 4
-8 l_ 1 l_ 4
+4 l_ 5
-4 l_ 1 l_ 5
+2 l_ 2 l_ 5
+2 l_ 4 l_ 5
\nonumber\\&&-l_5^2
-p_1
+4 p_ 2
+2 p_ 3
-\frac{\pi ^2}{3}
 \\
\sum_{i=1}^{\infty } \frac{16^{-i} \binom{2 i}{i}}{(1+2 i)^2}&=& \frac{3 \sqrt{3} c_ 1}{2}-\frac{\pi ^2}{3 \sqrt{3}}+\frac{3 \sqrt{3}}{2}-1 \\
\sum_{i=1}^{\infty } \frac{16^{-i} \binom{2 i}{i} 
\sum_{j=1}^i \frac{1}{j}}{1+2 i}&=& -5 \sqrt{3} c_ 1
+\frac{16 C}{3}
-\frac{2 \pi  l_ 1}{3}
-5 \sqrt{3}
+\frac{10 \pi ^2}{9 \sqrt{3}}
 \\
\sum_{i=1}^{\infty } \frac{16^{-i} \binom{2 i}{i} 
\sum_{j=1}^i \frac{1}{1+2 j}}{i}&=& \frac{l_ 2 l_ 3}{2}
+2 l_ 4
-2 l_ 1 l_ 4
+2 l_ 5
-2 l_ 1 l_ 5
+l_2 l_ 5
-\frac{p_ 1}
{2}
+2 p_ 2
\nonumber\\&&+2
-\frac{2 \pi }{3}
 \\
\sum_{i=1}^{\infty } \frac{16^{-i} \binom{2 i}{i} 
\sum_{j=1}^i \frac{1}{1+2 j}}{1+2 i}&=& \frac{\sqrt{3} c_ 1}{2}
+\frac{\pi  l_ 1}{3}
+\frac{\sqrt{3}}{2}
-\frac{\pi }{3}
-\frac{\pi ^2}{9 \sqrt{3}}
\end{eqnarray}
Weight 3:
\begin{eqnarray}
\sum_{i=1}^{\infty } \frac{1}{i^3 \binom{2 i}{i}}&=& \sqrt{3} \pi  c_ 1-\frac{4 \zeta_3}{3}-\frac{2 \pi ^3}{9 \sqrt{3}}+\sqrt{3} \pi  \label{w5rel3}\\
\sum_{i=1}^{\infty } \frac{
\sum_{j=1}^i \frac{1}{j}}{i^2 \binom{2 i}{i}}&=& \frac{\pi  c_ 1}{\sqrt{3}}-\frac{\zeta_3}{9}-\frac{2 \pi ^3}{27 \sqrt{3}}+\frac{\pi }{\sqrt{3}} \\
\sum_{i=1}^{\infty } \frac{
\sum_{j=1}^i \frac{1}{j^2}}{i \binom{2 i}{i}}&=& \sqrt{3} \pi  c_ 1-\frac{4 \zeta_3}{3}-\frac{35 \pi ^3}{162 \sqrt{3}}+\sqrt{3} \pi  \\
\sum_{i=1}^{\infty } \frac{2^i}{i^3 \binom{2 i}{i}}&=& \pi  C+\frac{\pi ^2 l_ 1}{8}-\frac{35 \zeta_3}{16} \label{w3rel1a}\\
\sum_{i=1}^{\infty } \frac{2^i 
\sum_{j=1}^i \frac{1}{j}}{i^2 \binom{2 i}{i}}&=& \frac{\pi ^2 l_ 1}{8}+\frac{7 \zeta_3}{16} \label{w3rel1b}\\
\sum_{i=1}^{\infty } \frac{2^i 
\sum_{j=1}^i \frac{1}{j^2}}{i \binom{2 i}{i}}&=& \pi  C+\frac{\pi ^2 l_ 1}{8}-\frac{35 \zeta_3}{16}+\frac{\pi ^3}{48} \\
\sum_{i=1}^{\infty } \frac{3^i}{i^3 \binom{2 i}{i}}&=& \frac{4 \pi  c_ 1}{\sqrt{3}}
+\frac{2 \pi ^2 l_ 2}{9}
-\frac{26 \zeta_3}{9}
+\frac{4 \pi }{\sqrt{3}}
-\frac{8 \pi ^3}
{27 \sqrt{3}}
 \\
\sum_{i=1}^{\infty } 
\frac{3^i 
\sum_{j=1}^i \frac{1}{j}}{i^2 \binom{2 i}{i}}&=& -\frac{2 \pi  c_ 1}{\sqrt{3}}
+\frac{2 \pi ^2 l_ 2}{9}
+\frac{13 \zeta_3}{9}
-\frac{2 \pi }{\sqrt{3}}
+\frac{4 \pi ^3}{27 \sqrt{3}}
 \\
\sum_{i=1}^{\infty } \frac{3^i 
\sum_{j=1}^i \frac{1}{j^2}}{i \binom{2 i}{i}}&=& \frac{4 \pi  c_ 1}{\sqrt{3}}
+\frac{2 \pi ^2 l_ 2}{9}
-\frac{26 \zeta_3}{9}
+\frac{4 \pi }{\sqrt{3}}
-\frac{4 \pi ^3}{27 \sqrt{3}}
 \\
\sum_{i=1}^{\infty } \frac{
\sum_{j=1}^i \frac{1}{j^2}}{(1+2 i) \binom{2 i}{i}}&=& \frac{\pi ^3}{81 \sqrt{3}} \\
\sum_{i=1}^{\infty } \frac{2^i 
\sum_{j=1}^i \frac{1}{j^2}}{(1+2 i) \binom{2 i}{i}}&=& \frac{\pi ^3}{48} \\
\sum_{i=1}^{\infty } \frac{3^i 
\sum_{j=1}^i \frac{1}{j^2}}{(1+2 i) \binom{2 i}{i}}&=& \frac{8 \pi ^3}{81 \sqrt{3}} \\
\sum_{i=1}^{\infty } \frac{
\sum_{j=1}^i \frac{1}{1+2 j}}{i^2 \binom{2 i}{i}}&=& -\frac{5 \pi  c_ 1}{2 \sqrt{3}}+\frac{35 \zeta_3}{18}+\frac{5 \pi ^3}{27 \sqrt{3}}-\frac{\pi }{2 \sqrt{3}}-4 \\
\sum_{i=1}^{\infty } \frac{2^i}{i^3 \binom{2 i}{i}}&=& \pi  C+\frac{\pi ^2 l_ 1}
{8}-
\frac{35 \zeta_3}{16} \\
\sum_{i=1}^{\infty } \frac{2^i 
\sum_{j=1}^i \frac{1}{1+2 j}}{i^2 \binom{2 i}{i}}&=& -\pi  C+\frac{7 \zeta_3}{2}-4+\pi  \\
\sum_{i=1}^{\infty } \frac{3^i}{i^3 \binom{2 i}{i}}&=& \frac{4 \pi  c_ 1}{\sqrt{3}}
+\frac{2 \pi ^2 l_ 2}{9}
-\frac{26 \zeta_3}{9}
+\frac{4 \pi }{\sqrt{3}}
-\frac{8 \pi ^3}{27 \sqrt{3}}
 \\
\sum_{i=1}^{\infty } \frac{3^i 
\sum_{j=1}^i \frac{1}{1+2 j}}{i^2 \binom{2 i}{i}}&=& -\frac{5 \pi  c_ 1}{\sqrt{3}}+\frac{91 \zeta_3}{18}+\frac{10 \pi ^3}{27 \sqrt{3}}-\frac{11 \pi }{3 \sqrt{3}}-4 \\
\sum_{i=1}^{\infty } \frac{16^{-i} \binom{2 i}{i}}{(1+2 i)^3}&=& \frac{7 \pi ^3}{216}-1 \\
\sum_{i=1}^{\infty } \frac{16^{-i} \binom{2 i}{i} 
\sum_{j=1}^i \frac{1}{1+2 j}}{(1+2 i)^2}&=& \frac{3}{8} \sqrt{3} c_ 1 l_ 1
+\frac{3 \sqrt{3} c_ 1}{4}
-\frac{5 \pi  c_ 1}{8}
+\frac{\sqrt{3} c_ 2}{2}
-\frac{c_ 3}{\sqrt{3}}
+\frac{5 c_ 4}{4 \sqrt{3}}
\nonumber\\&&-\frac{3}{4} \sqrt{3} c_ 6
+\frac{11 \sqrt{3} l_ 1}{8}
-\frac{7 \pi ^2 l_ 1}{36 \sqrt{3}}
-\frac{9}{8} \sqrt{3} l_ 2
+\frac{13 \zeta_3}
{24 \sqrt{3}}
\nonumber\\&&-\frac{3 \sqrt{3}}{2}
-\frac{\pi ^2}{12 \sqrt{3}}
+\frac{17 \pi ^3}{162}
 \\
\sum_{i=1}^{\infty } \frac{16^{-i} \binom{2 i}{i} 
\sum_{j=1}^i \frac{1}{(1+2 j)^2}}{1+2 i}&=& \frac{11 \pi ^3}{324}-\frac{\pi }{3}
\end{eqnarray}
Weight 4:
\begin{eqnarray}
 \sum_{i=1}^{\infty } 
\frac{1}{i^4 \binom{2 i}{i}}&=& \frac{17 \pi ^4}{3240} \\
\sum_{i=1}^{\infty } \frac{
\sum_{j=1}^i \frac{1}{j}}{i^3 \binom{2 i}{i}}&=& -\frac{3}{2} \sqrt{3} \pi  c_ 1 l_ 1
+18 c_ 1
+3 \sqrt{3} \pi  c_ 1
-\frac{61}{120} \pi ^2 c_ 1
-\frac{9 c_ 1^2}{2}
+\frac{2 \pi  c_ 2}{\sqrt{3}}
\nonumber\\&&-\frac{4 \pi  c_ 3}{3 \sqrt{3}}
+\frac{5 \pi  c_ 4}{3 \sqrt{3}}
-\sqrt{3} \pi  c_ 6
+\frac{1053 c_ 8}{80}
-2 c_ 9
-\frac{17 c_{10}}{15}
\nonumber\\&&-\frac{91 l_ 1 \zeta_3}{60}
+\frac{18 l_ 1}{5}
-\frac{\pi  l_ 1}{2 \sqrt{3}}
+\frac{5 \pi ^3 l_ 1}{27 \sqrt{3}}
+\frac{14}{135} \pi ^2 l_ 1^2
-\frac{14 l_ 1^4}{135}
\nonumber\\&&-\frac{3}{2} \sqrt{3} \pi  l_ 2
-\frac{112 p_ 4}{45}
-\frac{119 \zeta_3}{60}
-\frac{13 \pi  \zeta_3}{90 \sqrt{3}}
+\frac{981}{80}
+\frac{3 \sqrt{3} \pi }{5}
\nonumber\\&&-\frac{49 \pi ^2}{40}
-\frac{5 \pi ^3}{9 \sqrt{3}}
+\frac{1271 \pi ^4}{12150}
 \\
\sum_{i=1}^{\infty } \frac{
\sum_{j=1}^i \frac{1}{j^2}}{i^2 \binom{2 i}{i}}&=& \frac{7 \pi ^4}{1215} \\
\sum_{i=1}^{\infty } \frac{
\sum_{j=1}^i \frac{1}{j^3}}{i \binom{2 i}{i}}&=& -\frac{1}{6} \pi ^2 c_ 1
+\frac{27 c_ 8}
{4}
-
\frac{10 \pi  \zeta_3}{9 \sqrt{3}}
+\frac{27}{4}
-\frac{\pi ^2}{6}
-\frac{7 \pi ^4}{360}
 \\
\sum_{i=1}^{\infty } \frac{2^i}{i^4 \binom{2 i}{i}}&=& \pi  c_ 5
+2 C \pi  l_ 1
-\frac{1}{24} \pi ^2 l_ 1^2
+\frac{5 l_ 1^4}{48}
+\frac{5 p_ 4}{2}
-\frac{13 \pi ^4}{288}
 \\
 \sum_{i=1}^{\infty } \frac{2^i 
 \sum_{j=1}^i \frac{1}{j}}{i^3 \binom{2 i}{i}}&=& \pi  c_ 5
+C \pi  l_ 1
+2 C^2
+\frac{1}{12} \pi ^2 l_ 1^2
-\frac{l_ 1^4}{48}
-\frac{p_ 4}{2}
-\frac{11 \pi ^4}{360}
 \\
\sum_{i=1}^{\infty } \frac{2^i 
\sum_{j=1}^i \frac{1}{j^2}}{i^2 \binom{2 i}{i}}&=& \pi  c_ 5
+2 C \pi  l_ 1
-\frac{1}{24} \pi ^2 l_ 1^2
+\frac{5 l_ 1^4}{48}
+\frac{5 p_ 4}{2}
-\frac{49 \pi ^4}{1152}
 \\ 
 \sum_{i=1}^{\infty } \frac{2^i 
\sum_{j=1}^i \frac{1}{j^3}}{i \binom{2 i}{i}}&=& \pi  c_ 5
+6 c_ 7
+2 C \pi  l_ 1
-\frac{C \pi ^2}{4}
-\frac{1}{24} \pi ^2 l_ 1^2
+\frac{5 l_ 1^4}{48}
+\frac{5 p_ 4}{2}
\nonumber\\&&-\frac{13 \pi  \zeta_3}{16}
+6
-\frac{13 \pi ^4}{288}
 \\
\sum_{i=1}
^{\infty } 
\frac{3^i}{i^4 \binom{2 i}{i}}&=& \frac{7 \pi  c_ 1 l_ 1}{3 \sqrt{3}}
+\frac{4 \pi  c_ 1 l_ 2}{\sqrt{3}}
-\frac{10 \pi  c_ 1}{\sqrt{3}}
+\frac{31 \pi ^2 c_ 1}{45}
-\frac{4 \pi  c_ 2}{\sqrt{3}}
+\frac{8 \pi  c_ 3}{9 \sqrt{3}}
\nonumber\\&&-\frac{10 \pi  c_ 4}{9 \sqrt{3}}
-\frac{2 \pi  c_ 6}{\sqrt{3}}
-\frac{78 c_ 8}{5}
-\frac{64 c_{10}}{45}
+\frac{4 c_{11}}{9}
+\frac{1456 l_ 1 \zeta_3}{135}
\nonumber\\&&-\frac{26 l_ 2 \zeta_3}{9}
-\frac{128 l_ 1}{5}
-\frac{\pi  l_ 1}{3 \sqrt{3}}
-\frac{2 \pi ^3 l_ 1}{81 \sqrt{3}}
-\frac{64}{135} \pi ^2 l_ 1^2
+\frac{64 l_ 1^4}{135}
\nonumber\\&&+6 l_ 2
+\sqrt{3} \pi  l_ 2
-\frac{8 \pi ^3 l_ 2}{27 \sqrt{3}}
+\frac{1}{9} \pi ^2 l_ 2^2
+\frac{512 p_ 4}{45}
-\frac{112 \zeta_3}{45}
\nonumber\\&&+\frac{247 \pi  \zeta_3}{135 \sqrt{3}}
+\frac{54}{5}
+\frac{2 \sqrt{3} \pi }{5}
+\frac{4 \pi ^2}{15}
+\frac{2 \pi ^3}{27 \sqrt{3}}
-\frac{2449 \pi ^4}{18225}
 \label{w4rel3b}\\
\sum_{i=1}^{\infty } \frac{3^i 
\sum_{j=1}^i \frac{1}{j}}{i^3 \binom{2 i}{i}}&=& \frac{34 \pi  c_ 1 \
l_ 1}{3 \sqrt{3}}
-\frac{2 \pi  c_ 1 l_ 2}{\sqrt{3}}
+18 c_ 1
-\frac{28 \pi  c_ 1}{\sqrt{3}
}
+
\frac{179 \pi ^2 c_ 1}{360}
+\frac{45 c_ 1^2}{2}
\nonumber\\&&-\frac{8 \pi  c_ 2}{\sqrt{3}}
+\frac{32 \pi  c_ 3}{9 \sqrt{3}}
-\frac{40 \pi  c_ 4}{9 \sqrt{3}}
+\frac{4 \pi  c_ 6}{\sqrt{3}}
-\frac{429 c_ 8}{80}
+2 c_ 9
\nonumber\\&&+\frac{83 c_{10}}{45}
-\frac{2 c_{11}}{9}
-\frac{2093 l_ 1 \zeta_3}{540}
+\frac{13 l_ 2 \zeta_3}{9}
+\frac{46 l_ 1}{5}
+\frac{2 \pi  l_ 1}{3 \sqrt{3}}
\nonumber\\&&-\frac{32 \pi ^3 l_ 1}{81 \sqrt{3}}
+\frac{2}{15} \pi ^2 l_ 1^2
-\frac{2 l_ 1^4}{15}
-3 l_ 2
+2 \sqrt{3} \pi  l_ 2
+\frac{4 \pi ^3 l_ 2}{27 \sqrt{3}}
\nonumber\\&&+\frac{1}{9} \pi ^2 l_ 2^2
-\frac{16 p_ 4}{5}
+\frac{581 \zeta_3}{180}
-\frac{884 \pi  \zeta_3}{135 \sqrt{3}}
+\frac{27}{80}
+\frac{26 \sqrt{3} \pi }{5}
\nonumber\\&&-\frac{109 \pi ^2}{120}
+\frac{32 \pi ^3}{27 \sqrt{3}}
-\frac{11399 \pi ^4}{36450}
 \\
\sum_{i=1}^{\infty } \frac{3^i 
\sum_{j=1}^i \frac{1}{j^2}}{i^2 \binom{2 i}{i}}&=& \frac{7 \pi  c_ 1 \
l_ 1}{3 \sqrt{3}}
+\frac{4 \pi  c_ 1 l_ 2}{\sqrt{3}}
-\frac{10 \pi  c_ 1}{\sqrt{3}}
+\frac{31 \pi ^2 c_ 1}{45}
-\frac{4 \pi  c_ 2}{\sqrt{3}}
+\frac{8 \pi  c_ 3}{9 \sqrt{3}}
\nonumber\\&&-\frac{10 \pi  c_ 4}
{9 \sqrt{3}}
-\frac{2 \pi  c_ 6}{\sqrt{3}}
-\frac{78 c_ 8}{5}
-\frac{64 c_{10}}{45}
+\frac{4 c_{11}}{9}
+\frac{1456 l_ 1 \zeta_3}{135}
\nonumber\\&&-\frac{26 l_ 2 \zeta_3}{9}
-\frac{128 l_ 1}{5}
-\frac{\pi  l_ 1}{3 \sqrt{3}}
-\frac{2 \pi ^3 l_ 1}{81 \sqrt{3}}
-\frac{64}{135} \pi ^2 l_ 1^2
+\frac{64 l_ 1^4}{135}
\nonumber\\&&+6 l_ 2
+\sqrt{3} \pi  l_ 2
-\frac{8 \pi ^3 l_ 2}{27 \sqrt{3}}
+\frac{1}{9} \pi ^2 l_ 2^2
+\frac{512 p_ 4}{45}
-\frac{112 \zeta_3}{45}
\nonumber\\&&+\frac{247 \pi  \zeta_3}{135 \sqrt{3}}
+\frac{54}{5}
+\frac{2 \sqrt{3} \pi }{5}
+\frac{4 \pi ^2}{15}
+\frac{2 \pi ^3}{27 \sqrt{3}}
-\frac{2299 \pi ^4}{18225}
\label{w4rel2a} \\
\sum_{i=1}^{\infty } \frac{3^i 
\sum_{j=1}^i \frac{1}{j^3}}{i \binom{2 i}{i}}&=& \frac{7 \pi  c_ 1 l_ 1}{3 \sqrt{3}}
+\frac{4 \pi  c_ 1 l_ 2}{\sqrt{3}}
-\frac{10 \pi  c_ 1}{\sqrt{3}}
-\frac{29}{45} \pi ^2 c_ 1
-\frac{4 \pi  c_ 2}{\sqrt{3}}
+\frac{8 \pi  c_ 3}{9 \sqrt{3}}
\nonumber\\&&-\frac{10 \pi  c_ 4}{9 \sqrt{3}}
-\frac{2 \pi  c_ 6}{\sqrt{3}}
+\frac{12 c_ 8}{5}
-\frac{64 c_{10}}{45}
+\frac{4 c_{11}}{9}
+\frac{1456 l_ 1 \zeta_3}{135}
\nonumber\\&&-\frac{26 l_ 2 \zeta_3}{9}
-\frac{128 l_ 1}{5}
-\frac{\pi  l_ 1}{3 \sqrt{3}}
-\frac{2 \pi ^3 l_ 1}{81 \sqrt{3}}
-\frac{64}{135} \pi ^2 l_ 1^2
+\frac{64 l_ 1^4}{135}
\nonumber\\&&+6 l_ 2
+\sqrt{3} \pi  l_ 2
-\frac{8 \pi ^3 l_ 2}{27 \sqrt{3}}
+\frac{1}{9} \pi ^2 l_ 2^2
+\frac{512 p_ 4}{45}
-\frac{112 \zeta_3}{45}
\nonumber\\&&+\frac{187 \pi  \zeta_3}{135 \sqrt{3}}
+\frac{144}{5}
+\frac{2 \sqrt{3} \pi }{5}
-\frac{16 \pi ^2}{15}
+\frac{2 \pi ^3}{27 \sqrt{3}}
-\frac{2449 \pi ^4}{18225}
 \label{w4rel2b}\\
\sum_{i=1}^{\infty } \frac{
\sum_{j=1}^i \frac{1}{j^3}}{(1+2 i) \binom{2 i}{i}}&=& -\frac{1}{3} \pi ^2 c_ 1
+\frac{27 c_ 8}{2}
-\frac{20 \pi  \zeta_3}{9 \sqrt{3}}
+\frac{27}{2}
-\frac{\pi ^2}{3}
-\frac{4 \pi ^4}{81}
 \\
\sum_{i=1}^{\infty } \frac{2^i 
\sum_{j=1}^i \frac{1}{j^3}}{(1+2 i) \binom{2 i}{i}}&=& 6 c_ 7-\frac{\pi ^2 C}{4}-\frac{13 \pi  \zeta_3}{16}+6 \\
\sum_{i=1}^{\infty } \frac{3^i 
\sum_{j=1}^i \frac{1}{j^3}}{(1+2 i) \binom{2 i}{i}}&=& -\frac{8}{9} \pi ^2 c_ 1
+12 c_ 8
-\frac{8 \pi  \zeta_3}{27 \sqrt{3}}
+12
-\frac{8 \pi ^2}{9}
 \\
\sum_{i=1}^{\infty } \frac{
\sum_{j=1}
^i 
\frac{1}{1+2 j}}{i^3 \binom{2 i}{i}}&=& -\frac{3}{4} \sqrt{3} \pi  c_ 1 l_ 1
+\frac{45 c_ 1}{2}
+\frac{3}{2} \sqrt{3} \pi  c_ 1
-\frac{301}{240} \pi ^2 c_ 1
+\frac{9 c_ 1^2}{2}
+\frac{\pi  c_ 2}{\sqrt{3}}
\nonumber\\&&-\frac{2 \pi  c_ 3}{3 \sqrt{3}}
+\frac{5 \pi  c_ 4}{6 \sqrt{3}}
-\frac{1}{2} \sqrt{3} \pi  c_ 6
+\frac{1053 c_ 8}{160}
-c_9
-\frac{17 c_{10}}{30}
\nonumber\\&&-\frac{91 l_ 1 \zeta_3}{120}
+\frac{9 l_ 1}{5}
-\frac{\pi  l_ 1}{4 \sqrt{3}}
+\frac{5 \pi ^3 l_ 1}{54 \sqrt{3}}
+\frac{7}{135} \pi ^2 l_ 1^2
-\frac{7 l_ 1^4}{135}
\nonumber\\&&-\frac{3}{4} \sqrt{3} \pi  l_ 2
-\frac{56 p_ 4}{45}
-\frac{119 \zeta_3}{120}
-\frac{13 \pi  \zeta_3}{180 \sqrt{3}}
+\frac{3341}{160}
-\frac{31 \pi }{10 \sqrt{3}}
\nonumber\\&&-\frac{1241 \pi ^2}{720}
-\frac{5 \pi ^3}{18 \sqrt{3}}
+\frac{7919 \pi ^4}{97200}
 \\
\sum_{i=1}^{\infty } \frac{2^i 
\sum_{j=1}^i \frac{1}{1+2 j}}{i^3 \binom{2 i}{i}}&=& -\pi  c_ 5
-2 C \pi  l_ 1
+2 C^2
+\frac{1}{6} \pi ^2 l_ 1^2
-\frac{l_ 1^4}{6}
-4 p_ 4
+8
\nonumber\\&&-2 \pi 
-\frac{\pi ^2}
{4}
+
\frac{151 \pi ^4}{2880}
 \\
\sum_{i=1}^{\infty } \frac{3^i}{i^4 \binom{2 i}{i}}&=& \frac{7 \pi  c_ 1 l_ 1}{3 \sqrt{3}}
+\frac{4 \pi  c_ 1 l_ 2}{\sqrt{3}}
-\frac{10 \pi  c_ 1}{\sqrt{3}}
+\frac{31 \pi ^2 c_ 1}{45}
-\frac{4 \pi  c_ 2}{\sqrt{3}}
+\frac{8 \pi  c_ 3}{9 \sqrt{3}}
\nonumber\\&&-\frac{10 \pi  c_ 4}{9 \sqrt{3}}
-\frac{2 \pi  c_ 6}{\sqrt{3}}
-\frac{78 c_ 8}{5}
-\frac{64 c_{10}}{45}
+\frac{4 c_{11}}{9}
+\frac{1456 l_ 1 \zeta_3}{135}
\nonumber\\&&-\frac{26 l_ 2 \zeta_3}{9}
-\frac{128 l_ 1}{5}
-\frac{\pi  l_ 1}{3 \sqrt{3}}
-\frac{2 \pi ^3 l_ 1}{81 \sqrt{3}}
-\frac{64}{135} \pi ^2 l_ 1^2
+\frac{64 l_ 1^4}{135}
\nonumber\\&&+6 l_ 2
+\sqrt{3} \pi  l_ 2
-\frac{8 \pi ^3 l_ 2}{27 \sqrt{3}}
+\frac{1}{9} \pi ^2 l_ 2^2
+\frac{512 p_ 4}{45}
-\frac{112 \zeta_3}{45}
\nonumber\\&&+\frac{247 \pi  \zeta_3}{135 \sqrt{3}}
+\frac{54}{5}
+\frac{2 \sqrt{3} \pi }{5}
+\frac{4 \pi ^2}{15}
+\frac{2 \pi ^3}{27 \sqrt{3}}
-\frac{2449 \pi ^4}{18225}
 \\
\sum_{i=1}^{\infty } \frac{3^i 
\sum_{j=1}
^i 
\frac{1}{1+2 j}}{i^3 \binom{2 i}{i}}&=& \frac{13 \pi  c_ 1 l_ 1}{6 \sqrt{3}}
-\frac{5 \pi  c_ 1 l_ 2}{\sqrt{3}}
+15 c_ 1
+\frac{\pi  c_ 1}{\sqrt{3}}
-\frac{59}{48} \pi ^2 c_ 1
+\frac{57 c_ 1^2}{4}
\nonumber\\&&+\frac{2 \pi  c_ 2}{\sqrt{3}}
+\frac{4 \pi  c_ 3}{9 \sqrt{3}}
-\frac{5 \pi  c_ 4}{9 \sqrt{3}}
+\frac{5 \pi  c_ 6}{\sqrt{3}}
+\frac{663 c_ 8}{32}
+c_9
+\frac{55 c_{10}}{18}
\nonumber\\&&-\frac{7 c_{11}}{9}
-\frac{3913 l_ 1 \zeta_3}{216}
+\frac{91 l_ 2 \zeta_3}{18}
+43 l_ 1
+\frac{5 \pi  l_ 1}{6 \sqrt{3}}
-\frac{13 \pi ^3 l_ 1}{81 \sqrt{3}}
\nonumber\\&&+\frac{7}{9} \pi ^2 l_ 1^2
-\frac{7 l_ 1^4}{9}
-\frac{21 l_ 2}{2}
+\frac{\pi  l_ 2}{2 \sqrt{3}}
+\frac{10 \pi ^3 l_ 2}{27 \sqrt{3}}
-\frac{56 p_ 4}{3}
\nonumber\\&&+\frac{385 \zeta_3}{72}
-\frac{325 \pi  \zeta_3}{54 \sqrt{3}}
-\frac{161}{32}
+\frac{10 \pi }{3 \sqrt{3}}
-\frac{251 \pi ^2}{144}
+\frac{13 \pi ^3}{27 \sqrt{3}}
\nonumber\\&&+\frac{899 \pi ^4}{14580}
 \\
\sum_{i=1}^{\infty } \frac{16^{-i} \binom{2 i}{i}}{(1+2 i)^4}&=& \frac{27 \sqrt{3} c_ 8}{32}+\frac{\pi  \zeta_3}{12}-\frac{\pi ^4}{72 \sqrt{3}}+\frac{27 \sqrt{3}
}{32}-1  \label{w4rel1a}\\
\sum_{i=1}^{\infty } \frac{16^{-i} \binom{2 i}{i} 
\sum_{j=1}^i \frac{1}{1+2 j}}{(1+2 i)^3}&=& \frac{3 \sqrt{3} c_ 8}{16}
+\frac{7 \pi ^3 l_ 1}{216}
+\frac{\pi  \zeta_3}{24}
+\frac{3 \sqrt{3}}{16}
-\frac{7 \pi ^3}{216}
-\frac{\pi ^4}{324 \sqrt{3}}
 \label{w4rel1b}\\
\sum_{i=1}^{\infty } \frac{16^{-i} \binom{2 i}{i} 
\sum_{j=1}^i \frac{1}{(1+2 j)^2}}{(1+2 i)^2}&=& -\frac{3}{2} \sqrt{3} c_ 1
+\frac{\pi ^2 c_ 1}{16 \sqrt{3}}
+\frac{9 \sqrt{3} c_ 8}{16}
+\frac{\pi  \zeta_3}{9}
-\frac{15 \sqrt{3}}{16}
\nonumber\\&&+\frac{19 \pi ^2}{48 \sqrt{3}}
-\frac{\pi ^4}{72 \sqrt{3}}
 \\
\sum_{i=1}^{\infty } \frac{16^{-i} \binom{2 i}{i} 
\sum_{j=1}^i \frac{1}{(1+2 j)^3}}{1+2 i}&=& \frac{5 \pi  \zeta_3}{18}-\frac{\pi }{3}
\end{eqnarray}
For weight 5 and 6 we just list a few nice base identities

\begin{eqnarray}
\sum_{i=1}^{\infty } \frac{1}{i^5 \binom{2 i}{i}}&=& \frac{9}{4} \sqrt{3} \pi  c_ 8
-\frac{19 \zeta_5}{3}
+\frac{\pi ^2 \zeta_3}{9}
+\frac{9 \sqrt{3} \pi }{4}
-\frac{\pi ^5}{27 \sqrt{3}}
 \label{w5rel4}\\
\sum_{i=1}^{\infty } \frac{
\sum_{j=1}^i \frac{1}{j}}{i^4 \binom{2 i}{i}}&=&\frac{1}{2} \sqrt{3} \pi  c_ 8
-\frac{37 \zeta_5}{27}
+\frac{\pi ^2 \zeta (3)}{18}
+\frac{\sqrt{3} \pi }{2}
-\frac{2 \pi ^5}{243 \sqrt{3}}  \label{w5rel5}\\
\sum_{i=1}^{\infty } \frac{
\sum_{j=1}^i \frac{1}{1+2 j}}{i^4 \binom{2 i}{i}}&=&  -2 \sqrt{3} \pi  c_ 1
-\frac{17}{4} \sqrt{3} \pi  c_ 8
+\frac{713 \zeta_5}{54}
+\frac{8 \zeta (3)}{3}
-\frac{7 \pi ^2 \zeta (3)}{36}
\nonumber\\&&-16-\frac{43 \pi }{4 \sqrt{3}}
+\frac{2 \pi ^2}{9}
+\frac{4 \pi ^3}{9 \sqrt{3}}
+\frac{17 \pi ^5}{243 \sqrt{3}}\label{w5rel6}\\
\sum_{i=1}^{\infty } \frac{
\sum_{j=1}^i \frac{1}{j^2}}{i^3 \binom{2 i}{i}}&=& \frac{\pi ^3 c_1}{18 \sqrt{3}}
+\frac{3}{2} \sqrt{3} \pi  c_ 8
-\frac{142 \zeta_5}{27}
+\frac{4 \pi ^2 \zeta_3}{27}
+\frac{3 \sqrt{3} \pi }{2}
+\frac{\pi ^3}{18 \sqrt{3}}
\nonumber\\&&-\frac{7 \pi ^5}{243 \sqrt{3}}
 \\
\sum_{i=1}^{\infty } \frac{
\sum_{j=1}^i \frac{1}{j^3}}{i^2 \binom{2 i}{i}}&=& \frac{\zeta_5}{9}
+\frac{\pi ^2 \zeta_3}{27} 
\end{eqnarray}

\begin{eqnarray}
\sum_{i=1}^{\infty } \frac{16^{-i} \binom{2 i}{i}}{(1+2 i)^6}&=& \frac{297 \sqrt{3} c_{12}}{512}+\frac{\pi  \zeta_5}{16}+\frac{7 \pi ^3 \zeta_3}{864}-1+\frac{297 \sqrt{3}}{512}\nonumber\\&&-\frac{143 \pi ^6}{155520 \sqrt{3}}\\
\sum_{i=1}^{\infty } \frac{16^{-i} \binom{2 i}{i} \sum_{j=1}^i \frac{1}{(1+2 j)^3}}{(1+2 i)^3}&=&  -\frac{9}{64} \sqrt{3} c_{12}-\frac{11 \pi  \zeta_5}{432}+\frac{7 \pi ^3 \zeta_3}{216}-\frac{9 \sqrt{3}}{64}-\frac{7 \pi ^3}{216}\nonumber\\&&+\frac{13 \pi ^6}{58320 \sqrt{3}}
\end{eqnarray}

In the following we state some identities that can be derived from these base identities.
Combining (\ref{w3rel1a}) and (\ref{w3rel1b}) leads for instance to:
\begin{eqnarray}
 \sum_{k=1}^{\infty } \frac{2^k(\S{1}{k}-1/k)}{k^2 \binom{2 k}{k}}&=&\frac{21 \zeta_3-8 \pi  C}{8} 
\end{eqnarray}
From (\ref{w4rel3b}) and (\ref{w4rel2a}) we get:
\begin{equation}
\sum_{i=1}^{\infty } \frac{3^i \big(
        -\frac{1}{i^2}
        +
        \sum_{j=1}^i \frac{1}{j^2}
\big)}{i^2 \binom{2 i}{i}}=\frac{2 \pi ^4}{243}
\end{equation}
From (\ref{w4rel1a}) and (\ref{w4rel1b}) we find:
\begin{equation}
\sum_{i=1}^{\infty } \frac{16^{-i} \binom{2 i}{i} \left(\sum_{j=1}^i \frac{1}{1+2 j}-\frac{2}{9 (2i +1)}\right)}{(1+2 i)^3}=\frac{1}{216} \left(7 \pi ^3 l_1+5 \pi  \zeta_3-7 \pi ^3+48\right)
\end{equation}
From (\ref{w4rel2a}) and (\ref{w4rel2b}) we discover:
\begin{equation}
\sum_{i=1}^{\infty } \frac{3^i \big(
        \sum_{j=1}^i \frac{1}{j^3}
        -\frac{
        \sum_{j=1}^i \frac{1}{j^2}}{i}
\big)}{i \binom{2 i}{i}}= -\frac{4}{3} \pi ^2 c_ 1
+18 c_ 8
-\frac{4 \pi  \zeta_3}{9 \sqrt{3}}
+18
-\frac{4 \pi ^2}{3}
-\frac{2 \pi ^4}{243}
\end{equation}
From (\ref{w5rel1}), (\ref{w5rel2}), (\ref{w5rel3}), (\ref{w5rel4}), (\ref{w5rel5}) and (\ref{w5rel6}) we produce:
\begin{equation}
\sum_{i=1}^{\infty } \frac{-24
+\frac{2}{i^4}
+\frac{2}{i^2}
-\frac{4}{i}
-\frac{
\sum_{j=1}^i \frac{1}{j}}{2 i^3}
+\frac{
\sum_{j=1}^i \frac{1}{1+2 j}}{i^3}
}{i \binom{2 i}{i}}=\frac{11 \zeta_5}{9}-16
\end{equation}

In addition, we were able to prove the following conjectures from \cite{ZhiWei:2014}:

\begin{eqnarray}
\sum_{k=0}^\infty\f{\bi{2k}k}{(2k+1)16^k}\l(3\S{1}{2k+1}+\f4{2k+1}\r)&=&8C.
\end{eqnarray}

\begin{eqnarray}
\sum_{k=1}^\infty\f{\S{1}{2k}+2/(3k)}{k^2\bi{2k}k}&=&\zeta_3,\\
\sum_{k=1}^\infty\f{\S{1}{2k}+2\S{1}{k}}{k^2\bi{2k}k}&=&\f53\zeta_3,\\
\sum_{k=1}^\infty\f{\S{1}{2k}+17\S{1}{k}}{k^2\bi{2k}k}&=&\f52\sqrt3\ \pi K,\\
\sum_{k=1}^\infty\f{2^k}{k^2\bi{2k}k}\l(2\S{1}{2k}-3\S{1}{k}+\f2k\r)&=&\f 74\zeta_3,\\
\sum_{k=1}^\infty\f{2^k}{k^2\bi{2k}k}\l(6\S{1}{2k}-11\S{1}{k}+\f 8k\r)&=&2\pi C,\\
\sum_{k=1}^\infty\f{2^k}{k^2\bi{2k}k}\l(2\S{1}{2k}-7\S{1}{k}+\f2k\r)&=&-\f{\pi^2}2\log2,\\
\sum_{k=1}^\infty\f{3^k}{k^2\bi{2k}k}\l(6\S{1}{2k}-8\S{1}{k}+\f5k\r)&=&\f{26}3\zeta_3,\\
\sum_{k=1}^\infty\f{3^k}{k^2\bi{2k}k}\l(6\S{1}{2k}-10\S{1}{k}+\f7k\r)&=&2\sqrt3\,\pi K,\\
\sum_{k=1}^\infty\f{3^k}{k^2\bi{2k}k}\l(\S{1}{k}+\f1{2k}\r)&=&\f{\pi^2}3\log3.
\end{eqnarray}

\begin{eqnarray}
\sum_{k=0}^\infty\f{\bi{2k}k}{(2k+1)16^k}\sum_{j=0}^k\f1{(2j+1)^3}&=&\f 5{18}\pi\zeta_3.
\end{eqnarray}

\begin{eqnarray}
\sum_{k=1}^\infty\f{\S{1}{2k}-\S{1}{k}+2/k}{k^4\bi{2k}k}&=&\f{11}9\zeta_5,\\
\sum_{k=1}^\infty\f{3\S{2}{2k}-102\S{1}{k}+28/k}{k^4\bi{2k}k}&=&-\f{55}18\pi^2\zeta_3,\\
\sum_{k=1}^\infty\f{97\S{2}{2k}-163\S{1}{k}+227/k}{k^4\bi{2k}k}&=&\f{165}8\sqrt{3}\pi L,\\
\sum_{k=1}^\infty\f{\S{3}{k}}{k^2\bi{2k}k}&=&\f{\zeta_5+2\zeta_2\zeta_3}9,\\
\sum_{k=1}^\infty\f{\bi{2k}k}{(2k+1)16^k}\left(3\sum_{j=0}^k\f{1}{(2j+1)^4}-\f{1}{(2k+1)^4}\right)&=&\f{121 \pi^5}{17280}.
\end{eqnarray}

\begin{eqnarray}
\sum_{k=1}^\infty\f{\S{2}{k-1}-1/k^2}{k^4\bi{2k}k}&=&-\f{313 \pi^6}{612360},\\
\sum_{k=1}^\infty\f{3\S{4}{k}-1/k^4}{k^2\bi{2k}k}&=&\f{163 \pi^6}{136080},\\
\sum_{k=0}^\infty\f{\bi{2k}k}{(2k+1)16^k}\left(\sum_{j=0}^k\f{33}{(2j+1)^5}+\f{4}{(2k+1)^5}\right)&=&\f{35}{288}\pi^3\zeta_3+\f{1003}{96}\pi\zeta_5,\\
\sum_{k=0}^\infty\f{\bi{2k}k}{(2k+1)^3 16^k}\left(\sum_{j=0}^k\f{33}{(2j+1)^3}+\f{8}{(2k+1)^3}\right)&=&\f{245}{216}\pi^3\zeta_3-\f{49}{144}\pi\zeta_5.
\end{eqnarray}

Note that here $L:=2 c_8-\frac{8 \pi ^4}{729}+2$ and $K:=2 c_1-\frac{4 \pi ^2}{27}+2$.

\subsection*{Acknowledgements}
I would like to thank C. Schneider and F. Chyzak for calling my attention to \cite{ZhiWei:2014}. Additionally, I want to thank C. Schneider for useful discussions.



\end{document}